\documentclass[10pt]{amsart}
\usepackage{amssymb}
\usepackage{mathtools} % to get \coloneqq

\usepackage{geometry}                % See geometry.pdf to learn the layout options. There are lots.
\usepackage{graphicx}
\usepackage{epstopdf}
\usepackage{hyperref}

\usepackage{amsrefs}
\BibSpec{webpage}{%
  +{}{\PrintAuthors} {author}
  +{,}{ \textit} {title}
  +{}{ \parenthesize} {date}
  +{,}{ \url} {address}
  +{,}{ } {note}
  +{.}{ } {transition}
  +{}{ \parenthesize} {accessdate}
}

\newcommand{\NN}{\ensuremath{\mathbb N}}

\newcommand{\QQ}{\ensuremath{\mathbb Q}}

\newcommand{\scA}{{\mathcal A}}

\newcommand{\scF}{{\mathcal F}}

\newcommand{\floor}[1]{\lfloor {#1} \rfloor}

\newcommand{\dfloor}[1]{ \left\lfloor #1 \right \rfloor }

\newtheorem{theorem}{Theorem} %[section]

\newtheorem{corollary}[theorem]{Corollary}
\theoremstyle{definition}

\begin{document}

\title{Sets of natural numbers with proscribed subsets}

\author[O'Bryant]{Kevin O'Bryant}
\address{Department of Mathematics\\
College of Staten Island (CUNY)\\
Staten Island, NY 10314}
\thanks{Support for this project was provided by a PSC-CUNY Award,
  jointly funded by The Professional Staff Congress and The City
  University of New York.}
\thanks{This work will appear in the Journal of Integer Sequences.}
\email{kevin@member.ams.org}

\subjclass[2010]{11B05 11B25, 11B75, 11B83, 05D10.}
\keywords{Geometric progression-free sequences, Ramsey theory.}

\date{\today}

\begin{abstract}
Fix $\mathcal{A} \subseteq 2^{\mathbb{N}}$, and let $G_{\mathcal{A} }(n)$ be the maximum cardinality of a subset $X$ of $\{1,2,\dots,n\}$ with $2^X \cap \mathcal{A}  = \emptyset$. We consider the general problem of giving upper bounds on $G_{\mathcal{A} }(n)$, and give new results for some $\mathcal{A} $ that are closed under dilation. Specific examples addressed include sets that do not contain geometric progressions of length $k$ with integer ratio, sets that do not contain geometric progressions of length $k$ with rational ratio, and sets of integers that do not contain multiplicative squares, i.e., sets of the form $\{a,a r,a s,a r s\}$.
\end{abstract} 
\maketitle
%%%%%%%%%%%%%%%%%%%%%%%%%%%%%%%%%%%%%%%%%%%%%%%%%%        End Title

\section{Introduction}
Let $\scA$ be a collection of subsets of the natural numbers
($\NN\coloneqq \{1,2,\dots\}$ and $[n]\coloneqq \{1,2,\dots,n\}$), which we call the
{\em proscribed sets}, and let $S_{\scA}$ be the collection of sets of
natural numbers that do not have any subsets that are an element of $\scA$. Many of the most notorious problems in combinatorial number theory can be expressed as asking for properties of the elements of $S_{\scA}$. For example, let 
	\[ {\mathcal{AP}}_k \coloneqq  \big\{ \{a,a+d,a+2d,\dots,a+(k-1)d\} : a\in \NN, d \in \NN \big \}, \]
and $S_{{\mathcal{AP}}_k}$ is the collection of $k$-free sets of natural numbers. Let
	\[{\mathcal {IDON}} \coloneqq  \big\{ \{a,b,c,d\} : a<b\leq c < d, a+d=b+c \big\} ,	\]
and $S_{\mathcal{IDON}}$ is the collection of Sidon sets. Although there are exceptions, most of the problems of this sort that have been studied over the past 50 years concern affinely invariant $\scA$, i.e., if $A\in \scA$, then so is $d\ast A+t=\{da+t : a\in A\}$.

Recently, a number of works concerning sets that do not contain any $k$-term geometric progressions have appeared. Specifically, let
	\[ {\mathcal{GP}}_k \coloneqq  \big\{ a \ast \{1,r,r^2,\dots,r^{k-1}\} : a \in \NN, 1\neq r \in \NN \big\}	\]
and 
	\[ \widehat{\mathcal{GP}}_k \coloneqq  \big\{ a\ast \{1,r,r^2,\dots,r^{k-1}\} : a \in \NN, 1\neq r \in \QQ_+ \big\}.	\]
Further, for each set $X\subseteq\NN$ define
	\[ G_\scA(X) \coloneqq  \max\{ |A| : A \in S_{\scA}, A\subseteq X \}. \]
We will generate upper bounds on $G_{{\mathcal{GP}}_k}$ and $G_{\widehat{\mathcal{GP}}_k}$ in terms of other well-known Ramsey numbers, including Szemer\'{e}di numbers, density Hales-Jewett numbers and Moser numbers. Then, we turn our attention to some other forbidden sets.

\section{Statement of Main Result}
We call $\scF_0,\scF_1,\scF_2,\dots$ a {\em grading of $[n]$} if 
	\begin{enumerate}
		\item $\scF_0=\big\{ \{1\}, \{2\}, \dots, \{n\} \big \}$ is the collection of all singletons, \label{Cond:F_0}
		\item each $\scF_i$ is a collection of pairwise disjoint subsets of $[n]$,  \label{Cond:F_i packs}
		\item for each $i$, and for each $f_i\in\scF_{i},f_{i+1}\in\scF_{i+1}$, either $f_i\subseteq f_{i+1}$ for $f_i\cap f_{i+1}=\emptyset$, \label{Cond:F_i unrefines}
		\item for each $i$, and for each $f,g$ in $\scF_i$,  we have $G_{\scA}(f)=G_{\scA}(g)$.\label{Cond:Ramsey}
	\end{enumerate}
If additionally
	\begin{enumerate}			\setcounter{enumi}{4}
		\item for each $i$, each $f\in\scF_{i+1}$ is the disjoint union of exactly $k$ members of $\scF_{i}$,\label{Cond:expansion}
	\end{enumerate}
then we say that the grading has \emph{expansion $k$}. If instead
	\begin{enumerate}			\setcounter{enumi}{5}
		\item for each $i$ and each $f_{i+1}\in\scF_{i+1}$, there is $f_i\in \scF_i$ and $r$ distinct elements with $f_{i+1}=f_i\cup\{x_1,\dots,x_r\}$, and none of $x_1,\dots,x_r$ are contained in any member of any of $\scF_1,\dots,\scF_{i}$, \label{Cond:growth}
	\end{enumerate}
then we say that the grading has \emph{growth $r$}.
Whenever we set a particular grading, we assume that $\scF_0$ is what Condition (\ref{Cond:F_0}) requires, and that $\scF_d=\emptyset$ for any $d$ we don't expressly set.

\begin{theorem}\label{thm:main}
Let $\scA$ be a collection of proscribed sets, and let $\scF_0,\scF_1,\dots$ be a grading of $[n]$ with expansion $k\geq 2$, and let $R_i\geq G_{\scA}(f_i)$ for $f_i\in\scF_i$. Then
	\[
	\frac{G_{\scA}([n])}{n} \leq 1 - \sum_{i=1}^\infty (kR_{i-1} - R_{i}) \,\frac{|\scF_{i}|}{n}.
	\]
\end{theorem}

\begin{theorem}\label{thm:main2}
Let $\scA$ be a collection of proscribed sets, and let $\scF_0,\scF_1,\dots$ be a grading of $[n]$ with growth $r\geq 1$, and let $R_i\geq G_{\scA}(f_i)$ for $f_i\in\scF_i$. Then
	\[
	\frac{G_{\scA}([n])}{n} \leq 1 - \sum_{i=1}^\infty (r+R_{i-1} - R_{i}) \,\frac{|\scF_{i}|}{n}.
	\]
\end{theorem}

Observe that if a grading has expansion greater than $2$ or growth greater than $1$, then $\scF_d=\emptyset$ for sufficiently large $d$. Consequently, the infinite sums in the above theorems are, for each particular $n$, actually finite. Also, observe that the quantities $kR_{i-1}-R_i$ and $r+R_{i-1}-R_i$ are guaranteed to be nonnegative under the hypotheses of the theorems, so that the upper bounds in the above theorems are valid even if the infinite sums are truncated.

\section{Corollaries}

\subsection{Geometric progressions with prime-power ratio}
Let $p$ be a prime, and $k\geq 3,$ and let
	$$\scA = \big\{ a\ast \{1,p^{s}, p^{2s},\dots,p^{(k-1)s}\} : a\in \NN, s\in \NN \},$$
the geometric progressions of length $k$ whose ratio is a power of $p$. Let 
	$$\scF_i=\big\{b \ast \{1,p,\dots,p^{i}\} : (p,b)=1,\, 1\leq b \leq n/p^i\big\}.$$
This is a grading of $[n]$ with growth 1. As a geometric progression in $\{b,pb,\dots,p^i b\}$ is an arithmetic progression in the exponents $0,1,\dots,i$, we have $G_\scA(f) = r_k(i+1)$ for each $f\in \scF_i$, where $r_k(n)$ is the maximum size of a subset of $[n]$ that does not contain $k$-term arithmetic progressions. As $|\scF_i|= \frac{n}{p^i}\frac{\varphi(p)}{p}+O(1)$ for $1\leq i \leq \log_p n$, and $|\scF_i|=0$ for $i>\log_p n$, Theorem~\ref{thm:main2} now gives the bound
	\begin{align*}
	\frac{G_{\scA}([n])}{n} 
		&\leq 1 - \sum_{i=1}^{\floor{\log_p(n)}} (1+r_k(i-1)-r_k(i)) \frac{\frac{n}{p^i}\frac{\varphi(p)}{p}+O(1)}{n} \\
		&\lesssim 1- \left(1-\frac 1p\right) \sum_{i=1}^\infty \frac{1+r_k(i-1)-r_k(i)}{p^i}.
	\end{align*}
In~\cite{OBryant}, the author showed that this upper bound is asymptotically sharp (fixed $p$, with $n\to\infty$) and, perhaps surprisingly, is provably an irrational number.

We note that a set that avoids $k$-term arithmetic progressions cannot have $k$ consecutive elements, and so $r_k(n) \leq n - \floor{n/k}$, while by Szemer\'{e}di's Theorem, $r_k(n)=o(n)$. Therefore, there is a least $n$ with $r_k(n)<  n-\floor{k/n}$, and this value gives the improvement over ``easy'' in the above bound. We are not aware of any work explicitly aimed at finding this $n$, and some computations suggest that it depends on the {\em multiplicative} structure of $k$ and $k-1$.

\subsection{Three-term geometric progressions with friable integer ratio, McNew's method}
During preparation of this work, the author became aware of recent work of Nathan McNew [personal communication], a small portion of which fits into this framework. We give here just the facts with little justification, and leave the interested reader to seek out McNew's work.

Let $1=s_1<s_2<\dots$ be the natural numbers whose prime factors are among $2=p_1,p_2,\dots,p_d$ (the first $d$ prime numbers, with product $P_d$). Let $\scA^{(d)}$ be the set of 3-term geometric progressions whose ratio is among $s_2,s_3,\dots$. McNew takes
	\[
	\scF_i = \big\{ b \ast \{s_1,s_2,\dots,s_{i+1}\} : (b,P_d)=1, 1\leq b \leq n/s_{i+1}\big\},
	\]
which defines a grading with growth 1. As
	\[
	|\scF_i| \approx \frac{\varphi(P_d)}{P_d} \frac{n}{s_{i+1}},
	\]
Theorem~\ref{thm:main2} gives
	\[
	\frac{G_{\scA^{(d)}}([n])}{n} \lesssim 1 - \frac{\varphi(P_d)}{P_d} \sum_{i=1}^\infty \frac{1+R_{i-1}+R_i}{s_{i+1}},
	\]
where $R_i$ is the largest possible size of a subset of $\{1,s_2,s_3,\dots,s_{i+1}\}$ that does not contain any 3-term geometric progression. McNew shows further that this bound is asymptotically sharp, and that as $d\to\infty$ this bound approaches the answer to the problem in the next subsection (for $k=3$).

\subsection{Geometric progressions with integer ratio}

Consider $\scA={\mathcal{GP}}_k$, $k\geq 3$, and let
	$$\scF_1 = \big\{ a\ast \{1,2,2^2,\dots,2^{k-1}\} : 1\leq a \leq n/2^{k-1}, (a,2)=1\big\},$$
Then $R_0=1$, $R_1=G_{{\mathcal{GP}}_k}(\{1,2,2^2,\dots,2^{k-1}\}) = k-1$, $|\scF_1| = \floor{n/2^k+1/2}$, and Theorem~\ref{thm:main} gives
	\[
	\frac{G_{{\mathcal{GP}}_k}([n])}{n} \leq 1 - (k-(k-1))  \frac{\floor{n/2^k+1/2}}{n} \leq 1-2^{-k} + o(1).
	\]
This replicates the bound given in~\cite{Brown}.

Let 
	\[
	\scF_1 = \big\{ a2^{k(\ell-1)} \ast \{1,2,2^2,\dots,2^{k-1}\} : 1\leq \ell \leq \frac{1+\log_2 n}{k}, 1\leq a \leq n/2^{k\ell-1}, (a,2)=1\big\}.
	\]
As above, we have $R_0=1$, $R_1=k-1$, and 
	\[
	|\scF_1| = \sum_{\ell=1}^{(1+\log_2 n)/k} \dfloor{\frac{n}{2^{k\ell}} + \frac 12} = \frac{n}{2^k-1} + O(\log n),
	\] 
whence Theorem~\ref{thm:main} gives
	\[
	\frac{G_{{\mathcal{GP}}_k}([n])}{n} \leq 1 - \frac{1}{2^k-1} + o(1).
	\]
This replicates the bound given in~\cite{Riddell}, and rediscovered in \cite{Beiglbock}.
	
Now, we go further, providing the first example of the power of Theorem~\ref{thm:main} and giving our first original corollary. Denote the sequence of prime numbers as $2=p_1<p_2<p_3<\dots$, and let $P_d\coloneqq \prod_{i=1}^d p_i$ and
	\[a_d \coloneqq  \left\{ \prod_{i=1}^d p_i^{e_i} : 0\leq e_i <k\right\}.\]
We set up our grading with expansion $k$ as follows: $\scF_0=\big\{ \{b\} : 1\leq b \leq n \}$, and for $d\geq 1$
	\[\scF_d \coloneqq  \left\{ 2^{k(\ell-1)}b \ast a_{d} : 1\leq \ell \leq 1+\frac{\log_2 (n/P_d^{k-1})}k, 1\leq b \leq \frac{n}{P_d^{k-1}2^{k(\ell-1)}}, (b,P_d)=1 \right\} .\] 
Before continuing, we must establish that this is actually a grading. Condition (\ref{Cond:F_0}) is immediate. The smallest element of $2^{k(\ell-1)}b\ast a_d$ is $2^{k(\ell-1)} b\geq 1$ and the largest element is 
	\[
	2^{k(\ell-1)}b P_d^{k-1} \leq 2^{k(\ell-1)}\cdot \frac{n}{P_d^{k-1} 2^{k(\ell-1)}} P_d^{k-1} = n,
	\]
which establishes that every member of $\scF_d$ is a subset of $[n]$. Fix $d\geq 1$, and suppose that 
	\[
	x \in 2^{k(\ell_1-1)}b_1 \ast a_d \cap 2^{k(\ell_2-1)}b_2 \ast a_d.
	\]
Let $v_p(x)$ be the highest exponent of $p$ that divides $x$. As $x\in 2^{k(\ell_1-1)}b_1 \ast a_d$, we see that $k(\ell_1-1) \leq v_2(x) < k(\ell_1-1)+k-1$, and likewise $k(\ell_2-1) \leq v_2(x) < k(\ell_2-1)+k$. Therefore, $\ell_1=\ell_2$, and so we have
	\[
	\frac{x}{2^{k(\ell-1)}} = b_1 \prod_{i=1}^d p_i ^{e_i} = b_2 \prod_{i=1}^d p_i^{e'_i}.
	\]
As $(b,P_d)=1$, we observe that $b_1=b_2$, and then that $e_i=e'_i$ for each $i$. In other words, the members of $\scF_d$ are pairwise disjoint, confirming condition (\ref{Cond:F_i packs}). Conditions (\ref{Cond:F_i unrefines}) and (\ref{Cond:expansion}) are both settled by the observation that
	\[
	2^{k(\ell-1)} b \ast a_d = \bigcup_{i=0}^{k-1} 2^{k(\ell-1)} (b p_d^i) \ast a_{d-1},
	\]
and so each member of $\scF_{i+1}$ is the disjoint union of exactly $k$ members of $\scF_i$. Condition (\ref{Cond:Ramsey}) follows from the observation that ${\mathcal{GP}}_k$ is closed under dilation, and so
	\[
	G_{{\mathcal{GP}}_k}(2^{k(\ell-1)}b \ast a_{d}) = G_{{\mathcal{GP}}_k}(a_d).
	\]

Suppose that $x_1,\dots,x_k$ is a geometric progression in $a_d$, with integer ratio 
	$$r=\frac{x_2}{x_1}= \prod_{i=1}^d p_i^{r_i}$$
with at least one $r_i>0$, and all $r_i\geq 0$. As $x_k=x_1 r^{k-1}\in a_d$, we see that no $r_i >1$. In particular, the sequence of $d$-tuples 
	$$ \big(v_{p_1}(x_i), v_{p_2}(x_i), \dots, v_{p_d}(x_i) \big),\qquad (1\leq i \leq k)$$
has some coordinates fixed, while the others count in unison from 0 up to $k-1$. That is, they are a combinatorial line in $\{0,1,\dots,k-1\}^d$, and each combinatorial line in $\{0,1\dots,k-1\}^d$ is generated by a geometric progression in $a_d$. As in the Polymath 1 project~\cite{Polymath}, we denote the largest possible size of a subset of  $\{0,1\dots,k-1\}^d$ that does not contain a combinatorial line as $c_{d,k}$. To wit,
	\[
	G_{{\mathcal{GP}}_k}(a_d) = c_{d,k}
	\]

For each value of $\ell$, $\scF_d$ has one member for each $b$ between 1 and $\frac{n}{P_d^{k-1}2^{k(\ell-1)}}$ that is relatively prime to $P_d$. The proportion of numbers relatively prime to $P_d$ is $\varphi(P_d)/P_d$, and so there are
	$$
	\frac{n}{P_d^{k-1}2^{k(\ell-1)}} \, \frac{\varphi(P_d)}{P_d} + O(P_d) = n \frac{\varphi(P_d)}{P_d^k} \cdot \frac{1}{2^{k(\ell-1)}}+ O(1)
	$$
such values of $b$. Summing this over $\ell$, with $M\coloneqq 1+\log_2(n/P_d^{k-1})/k$, yields
	\begin{align*}
	|\scF_d| &= \sum_{\ell=1}^M \left( n \frac{\varphi(P_d)}{P_d^k} \cdot \frac{1}{2^{k(\ell-1)}}+ O(1) \right) \\
	&=n\, \frac{\varphi(P_d)}{P_d^{k}}\sum_{\ell=1}^M \frac{1}{2^{k(\ell-1)}} + O(\log n) \\
	&=n \, \frac{\varphi(P_d)}{P_d^{k}} \, \frac{2^k}{2^k-1} + O(\log n).
	\end{align*}

Theorem~\ref{thm:main} gives
	\[
	\frac{G_{{\mathcal{GP}}_k}([n])}{n} \leq 1 - \frac{2^k}{2^k-1} \sum_{d=1}^\infty ( k c_{d-1,k} - c_{d,k}) \frac{\varphi(P_d)}{P_d^{k}}+o(1).
	\]
By the Density Hales-Jewett Theorem infinitely many of the $kc_{d-1,k}-c_{d,k}$ are positive. This bound is superior to any in the literature (the previous best corresponds to taking only the first term of the sum), and so we state it explicitly.

\begin{corollary}
Let $G_{{\mathcal{GP}}_k}([n])$ be the largest possible size of a subset of $[n]$ that does not contain any $k$-term geometric progression with integer ratio. Then
	\[
	\limsup_{n \to \infty} \frac{G_{{\mathcal{GP}}_k}([n])}{n} 
		\leq 1 - \frac{2^k}{2^k-1} \sum_{d=1}^\infty ( k c_{d-1,k} - c_{d,k}) \frac{\varphi(P_d)}{P_d^{k}},
	\]
where $P_d$ is product of the first $d$ primes, $\varphi$ is Euler's phi function, and $c_{d,k}$ are the density Hales-Jewett numbers.
\end{corollary}	

The Density Hales-Jewett Theorem states that $c_{d,k} = o(k^d)$ (with $k$ fixed, $d\to\infty$), and the recent Polymath project~\cite{Polymath} found $c_{d,3}$ for $d\leq 6$:
  \[c_{0,3}=1,c_{1,3}=2, c_{2,3}=6, c_{3,3}=18, c_{4,3}=52, c_{5,3}=150, c_{6,3}=450, 1302\leq c_{7,3} \leq 1348,\] 
which is A156989 in the OEIS. Using these values, we get
	\[
	\limsup_{n \to \infty} \frac{G_{{\mathcal{GP}}_3}([n])}{n} 
		< 1 - \frac 87 \left(1\, \frac{\varphi(2)}{2^3}+ 2\, \frac{\varphi(2\cdot3\cdot5\cdot 7)}{(2\cdot 3 \cdot 5 \cdot 7)^3} 
	                    + 6 \, \frac{\varphi(2\cdot3\cdot5\cdot 7\cdot 11)}{(2\cdot 3 \cdot 5 \cdot 7\cdot 11)^3} \right) < 0.857131.
	\]
Assorted other useful values were also computed \cites{HigherDimDHJ,UpperAndLowerDHJ} for $k>3$:
  \[c_{4,4}=183, c_{5,4}\leq 732, c_{4,6}\leq 1079,\]
although these calculations were not subjected to the same scrutiny and should be considered less reliable. They lead to positive, albeit numerically miniscule, 	improvements on the previously known bounds for $\limsup_n \frac{G_{{\mathcal{GP}}_k}([n])}{n}$ for $k=4$ and $k=6$. By the Density Hales-Jewett Theorem, infinitely many of the $kc_{d-1,k}-c_{d,k}$ are positive, and so this gives an improvement for all $k$, even though we are unable assess the magnitude of the improvement for $k\not\in \{3,4,6\}$.

\subsection{Geometric progressions with rational ratio}
As ${\mathcal{GP}}_k \subseteq \widehat{\mathcal{GP}}_k$, we know that
	\[G_{\widehat{\mathcal{GP}}_k}(X) \leq G_{{\mathcal{GP}}_k}(X) \]
for any $X$, and so the bounds of the previous section apply here, too. We can do a bit better, however, because not every geometric progression with rational ratio in $a_d$ lands on a combinatorial line in $\{0,1,\dots,k-1\}^d$. The appropriate structure is called a geometric line: the $k$ distinct points $x_i \in \{0,1,\dots,k-1\}^d$ are a geometric line if the coordinates fall into three categories, one where the coordinate value never changes, one where the coordinate value counts up from 0 to $k-1$, and one (possibly empty) where the coordinate value counts down from $k-1$ to 0. The largest possible size of a subset of $\{0,1,\dots,k-1\}^d$ that contains no geometric lines is denoted $c'_{d,k}$ and were also studied in the Polymath~\cite{Polymath} project:
	\[c'_{0,3}=1,c'_{1,3}=2, c'_{2,3} = 6,c'_{3,3} = 16,c'_{4,3} = 43,c'_{5,3} = 124,c'_{6,3} = 353,\]
which calls this the sequence of Moser numbers.
\begin{corollary}
Let $G_{{\mathcal{GP}}_k}([n])$ be the largest possible size of a subset of $[n]$ that does not contain any $k$-term geometric progression with integer or rational ratio. Then
	\[
	\limsup_{n \to \infty} \frac{G_{{\mathcal{GP}}_k}([n])}{n} 
		\leq 1 - \frac{2^k}{2^k-1} \sum_{d=1}^\infty ( k c'_{d-1,k} - c'_{d,k}) \frac{\varphi(P_d)}{P_d^{k}},
	\]
where $P_d$ is product of the first $d$ primes, $\varphi$ is Euler's phi function, and $c'_{d,k}$ are the Moser numbers.
\end{corollary}	

With $k=3$, this improves the previous best bound of $6/7$ to
	\[
	\limsup_{n \to \infty} \frac{G_{\widehat{\mathcal{GP}}_3}([n])}{n} < \frac 67 - \frac{16755239936}{23695945898625} < 0.856436.
	\]
We don't know any nontrivial Moser numbers with $k>3$, although the Density Hales-Jewett Theorem implies that infinitely many of the $k c'_{d-1,k} - c'_{d,k}$ must be positive.

\subsection{Geometric squares}
A geometric square is a set of 4 natural numbers of the form $\{a,ar,as,ars: a,r,s \in \NN\}$; set 
	\[
	\scA=\big\{ a \ast \{1,r,s,rs\} : a,r,s \in \NN, 1<r<s \big\}.
	\]
We note that for each $n$, the set
	\[
	\left(n \ast (\tfrac{1}{6},1]\right) \cap \NN
	\]
does not contain any geometric squares, and is a subset of $[n]$ with cardinality $\sim (5/6) n$.

Let $2=p_1<p_2<\cdots$ be the prime numbers, $P_d\coloneqq  \prod_{i=1}^d p_i$, and set
	\[a_d \coloneqq  \left\{ \prod_{i=1}^d p_i^{e_i} : 0\leq e_i \leq 1 \right\}.\]
We now define the grading (leaving the details to the reader)
	\[\scF_d \coloneqq  \big\{ b 4^i \ast a_d : (b,P_d)=1, 1\leq b \leq n/(P_d \cdot4^i), 0\leq i \leq \log_4(n/P_d) \big\}.\]
This grading has expansion 2 and
	\begin{align*}
	|\scF_d| &= \sum_{i=0}^{\log_4(n/P_d)} \left( \frac{n}{P_d \cdot 4^i} \frac{\varphi(P_d)}{P_d} + O(1) \right) \\
		&= n \frac{\varphi(P_d)}{P_d^2} \frac 43  + O(\log n).
	\end{align*}

The Ramsey numbers $G_{\scA}(a_d)$, which we will denote as $c_{d,2,2}$, deserved to have been studied before, but we are not aware of any such computation. In particular, $c_{d,2,2}$ is the greatest number $c$ for which there is a family of $c$ subsets of $[d]$, no four of which have the relations $A_1\subset A_2 \subset A_4,A_1 \subset A_3 \subset A_4$. More generally, we would set $c_{d,s,2}$ would be the maximum cardinality of a family of subsets of $[d]$ that does not contain a sub-family that is lattice isomorphic to the power set of $[s]$. Sperner's Theorem is equivalent to the assertion that $c_{d,1,2}= \binom{d}{\floor{d/2}}$. Even more generally, $c_{d,s,k}$ is the maximum cardinality of a subset of $[k]^d$ that does not contain a combinatorial space with dimension $s$.

We report the following values:
	\[c_{0,2,2}=1, c_{1,2,2} = 2, c_{2,2,2} = 3, c_{3,2,2}=6, c_{4,2,2}=11, c_{5,2,2}=21.\]

Theorem~\ref{thm:main} now gives us a clean bound.
\begin{corollary}
Let $G_{\scA}([n])$ be the largest possible size of a subset of $[n]$ that does not contain any subset of the form $\{a,ar,as,ars\}$ with $a,r,s$ being natural numbers with $1<r<s$. Then
	\[
	\frac 56 \leq \limsup_{n \to \infty} \frac{G_{\scA}([n])}{n} \leq  \frac{3699337}{4002075}
		< 1 - \frac 43 \sum_{d=1}^\infty ( 2 c_{d-1,2,2} - c_{d,2,2}) \frac{\varphi(P_d)}{P_d^2},
	\]
where $P_d$ is product of the first $d$ primes, $\varphi$ is Euler's phi function, and $c_{d,2,2}$ are the generalized Sperner numbers defined above.
\end{corollary}	

We are unaware if there a subset of $\NN$ that avoids $\scA$ and has positive density.	

\section{Proof of Main Result}
We prove Theorem~\ref{thm:main}, and leave the similar, and slightly easier, proof of Theorem~\ref{thm:main2} to the reader. First, note that since the grading has expansion at least 2, every member of $\scF_d$ must have size at least $2^d$ and be subsets of $[n]$; thus $\scF_d=\emptyset$ for $d> \log_2 n$. In particular, the infinite sum in the statement of the theorem is actually a finite sum.

For each $b\in [n]$, let $\delta(b)$ be the largest $d$ such that there is $A\in \scF_d$ with $b\in A$, and let $A_b$ be the unique set with $b\in A_b \in \scF_{\delta(b)}$. The function $\delta(b)$ is well defined as $b\in \{b\} \in \scF_0$ (condition~\ref{Cond:F_0}), and $\scF_d=\emptyset$ for sufficiently large $d$, and $A_b$ is well-defined as $\scF_{\delta(b)}$ is a collection of disjoint subsets of $[n]$ (condition~\ref{Cond:F_i packs}). Moreover, observe that if $A_b,A_c$ are not identical, then they have no intersection. For otherwise, if $x\in A_b \cap A_c$, then by Condition (\ref{Cond:F_i packs}), $\delta(b)\not= \delta(c)$, say $\delta(b)<\delta(c)$. By the expansion property, $A_c$ is the disjoint union of $k$ members of $\scF_{\delta(c)-1}$, and by induction is the disjoint union of $k^{\delta(c)-\delta(b)-1}$ members of $\scF_{\delta(b)+1}$. One of those members must contain $x$, and so by condition~\ref{Cond:F_i unrefines}, must contain all of $A_b$. In particular, $b \in A_c$, whence $\delta(b)\geq \delta(c)$ and so $A_b=A_c$.

Thus, the family
	\[	{\mathcal P} \coloneqq \big\{ A_b : b\in [n] \big\} 	\]
is a partition of $[n]$. Let $\alpha_i$ be the number of members of ${\mathcal P}$ that are in $\scF_i$. For disjoint sets $X,Y$, we have
	\[	G_{\scA}(X \cup Y) \leq G_{\scA}(X) + G_{\scA}(Y).	\]
Applying this principle to the partition ${\mathcal P}$, we have
	\[	G_{\scA}([n])  \leq \sum_{A \in {\mathcal P}} G_{\scA}(A)  \leq  \sum_{i=0}^{\log_2 n} \alpha_i R_i,	\]
using condition~\ref{Cond:Ramsey}. By high-school algebra, we have
	\[ \sum_{i=0}^{\log_2 n} \alpha_i R_i = R_0 \sum_{i=0}^{\log_2 n} k^i \alpha_i- \sum_{i=1}^{\log_2 n} (k R_{i-1}-R_i)\sum_{j=i}^{\log_2 n} k^{j-i} \alpha_j. \]
Note now that 	
	\[ R_0 \leq 1, \qquad \sum_{i=0}^{\log_2 n} k^i \alpha_i = |\scF_0| = n, \qquad \sum_{j=i}^{\log_2 n} k^{j-i} \alpha_j = |\scF_i|, \]
and the Theorem is proved.

\section{Remaining Problems}

Aside from the obvious ``sharpen the given bounds'' problems, we single out 3 interesting problems.

Let $r_k(n)$ be the largest possible size of a subset of $\{1,2,\dots,n\}$ that does not have a subset that is a $k$-term arithmetic progression. For each $k$, what is the least $n$ with $r_k(n) < n-\floor{n/k}$?

Is there a subset of $\NN$ that has positive density and does not contain a subset of the form $\{a,ar,as,ars\}$, where $a,r,s$ are natural numbers and $r,s$ are both greater than $1$?

Is there a clean formula for $c_{d,2,2}$, the maximum size of a family of subsets of $[d]$ that does not contain 4 sets $A_1,A_2,A_3,A_4$ with $A_1 \subset A_2 \subset A_4,A_1\subset A_3 \subset A_4$?

%\clearpage

\end{document}